\documentclass[11pt]{article}

\usepackage{amsfonts,amsmath,amsthm,amssymb}
\usepackage{latexsym}

\newtheorem{theorem}{Theorem}[section]
\newtheorem{corollary}[theorem]{Corollary}
\newtheorem{lemma}[theorem]{Lemma}
\newtheorem{example}[theorem]{Example}

\newtheorem{proposition}[theorem]{Proposition}
\newtheorem{fact}[theorem]{Fact}

\theoremstyle{definition}
\newtheorem{claim}[theorem]{Claim}
\newtheorem{eee}[theorem]{Example}
\newtheorem{rrr}[theorem]{Remark}
\newtheorem{sss}[theorem]{Statement}
\newtheorem{defi}[theorem]{Definition}
\newtheorem{notation}[theorem]{Notation}
\newtheorem{qqq}[theorem]{Question}

\newcommand{\beql}[1]{ \begin{equation} \label{#1} }
\newcommand{\eeq}{\end{equation}}
\newcommand{\comment}[1]{}

\newcommand{\Abs}[1]{{\left|{#1}\right|}}

\newcommand{\Prob}[1]{{{{\mathbb Pr}}\left[{#1}\right]}}
\newcommand{\Set}[1]{{\left\{{#1}\right\}}}

\newcommand{\RR}{{\mathbb R}}
\newcommand{\CC}{{\mathbb C}}
\newcommand{\ZZ}{{\mathbb Z}}

\newcommand{\QQ}{{\mathbb Q}}

\newcommand{\ft}[1]{\widehat{#1}}

\newcommand{\supp}{{\rm supp\,}}

\newcommand{\eat}[1]{}

\newcounter{rem}
\newcounter{othm}
\def\theothm{\Alph{othm}}

\newcommand{\bt}{\begin{theorem}}
\newcommand{\bl}{\begin{lemma}}
\newcommand{\bc}{\begin{claim}}
\newcommand{\bex}{\begin{eee}}
\newcommand{\br}{\begin{rrr}}
\newcommand{\bs}{\begin{sss}}
\newcommand{\bd}{\begin{defi} \upshape}
\newcommand{\bn}{\begin{notation} \upshape}
\newcommand{\bq}{\begin{qqq}}
\newcommand{\bcor}{\begin{corollary}}

\newcommand{\bp}{\noindent\textbf{Proof. }}

\newcommand{\et}{\end{theorem}}
\newcommand{\el}{\end{lemma}}
\newcommand{\ec}{\end{claim}}
\newcommand{\eex}{\end{eee}}
\newcommand{\erem}{\end{rrr}}
\newcommand{\es}{\end{sss}}
\newcommand{\ed}{\end{defi}}
\newcommand{\en}{\end{notation}}
\newcommand{\eq}{\end{qqq}}
\newcommand{\ecor}{\end{corollary}}

\newcommand{\ep}{\hspace{\stretch{1}}$\square$\medskip}

\newcommand{\beq}{\begin{equation}}

\newcommand{\N}{\mathbb{N}}
\newcommand{\Z}{\mathbb{Z}}
\newcommand{\Zn}{\mathbb{Z}_n}
\newcommand{\Q}{\mathbb{Q}}
\newcommand{\R}{\mathbb{R}}
\newcommand{\C}{\mathbb{C}}

\newcommand{\al}{\alpha}
\newcommand{\be}{\beta}
\newcommand{\ga}{\gamma}

\newcommand{\zz}{\zeta_{2k}}

\newcommand{\su}{\subset}
\newcommand{\sm}{\setminus}
\newcommand{\FTE}{\widehat{\chi_E}}
\newcommand{\FTF}{\widehat{\chi_F}}
\newcommand{\ftf}{\widehat{f}}
\newcommand{\ftg}{\widehat{g}}
\newcommand{\md}{{\rm{mod\,}}}   
\newcommand{\sFTE}{\supp\FTE}

\newcommand{\sftf}{\supp\ftf}

\newcommand{\class}[1]{\langle #1 \rangle_n}
\newcommand{\e}[2]{e^{-2\pi i\frac{#1}{#2}}}
\newcommand{\tor}{\textrm{ or }}
\newcommand{\tand}{\textrm{ and }}
\newcommand{\tif}{{\ \rm{ if }\ }}
\newcommand{\Range}{{\rm{ Range }}}

\begin{document}

\title{On the determination of sets by their triple correlation in finite cyclic groups}

\author{
{\sc{Tam\'as Keleti}\thanks{
Department of Analysis, E\"otv\"os Lor\'and University,
P\'azm\'any P\'eter s\'et\'any 1/C, H-1117 Budapest, Hungary.
E-mail: {\tt elek@cs.elte.hu}.
Partially supported by 
OTKA grants 049786 and F 43620 and by
European Commission IHP Network HARP
(Harmonic Analysis and Related Problems), Contract Number: HPRN-CT-2001-00273 - HARP.
This research started when the first author visited
the second author at the University of Crete and it continued while
the first author was a visitor at the Alfr\'ed R\'enyi Institute of
Mathematics of the Hungarian Academy of Science.
}}\and{{\sc Mihail N. Kolountzakis}\thanks{Department of Mathematics, Univ.\ of Crete, GR-71409 Iraklio, Greece. E-mail: {\tt kolount@gmail.com}.
Partially supported by European Commission IHP Network HARP
(Harmonic Analysis and Related Problems), Contract Number: HPRN-CT-2001-00273 - HARP,
and by grant INTAS 03-51-5070 (2004) (Analytical and Combinatorial Methods in Number Theory and Geometry).
Also by the Greek research program "Pythagoras 2" (75\% European funds
and 25\% National funds).
}
}
}

\date{16 March 2006}

\maketitle

\begin{abstract}
Let $G$ be a finite abelian group and $E$ a subset of it.
Suppose that we know for all subsets $T$ of $G$
of size up to $k$ for how many $x \in G$ the translate $x+T$
is contained in $E$.
This information is collectively called the $k$-deck of $E$.
One can naturally extend the domain of definition of the $k$-deck to include
functions on $G$.
Given the group $G$ when is the $k$-deck of a set in $G$ sufficient to determine
the set up to translation?
The $2$-deck is not sufficient (even when we allow for
reflection of the set, which does not change the $2$-deck)
and the first interesting case is $k=3$.
We further restrict $G$ to be cyclic and determine the values of $n$ for which the
$3$-deck of a subset of $\ZZ_n$ is sufficient to determine the set up to translation.
This completes the work begun by Gr\"unbaum and Moore \cite{GM} as far as the $3$-deck
is concerned.
We additionally estimate from above the probability that for a random subset of $\ZZ_n$
there exists another subset, not a translate of the first, with the same $3$-deck.
We give an exponentially small upper bound when the previously known
one was $O(1\bigl / \sqrt{n})$.  
\end{abstract}

\section{Introduction to the problem and results}
\label{sec:introduction}

Let $G$ be a finite abelian group, written additively,
and $f:G \to \RR$ 
be a function.
For $k \ge 2$ we define the $k$-deck or $k$-th order correlation of $f$ as the function
$$
N_{f,k}:G^{k-1} \to \RR
$$
defined by
\begin{equation}\label{k-deck-defn}
N_{f,k}(x_1,\ldots,x_{k-1}) = \sum_{x \in G} f(x) f(x+x_1) \cdots f(x+x_{k-1}).
\end{equation}
When $E \subseteq G$ and $f(x) = \chi_E(x)$ is the indicator function of $E$
we also write $N_{E,k}$ in place of $N_{f,k}$.
In this case, of $f = \chi_E$, it is easy to see that the number
$N_{E,k}(x_1,x_2,\ldots,x_{k-1})$ is precisely the number of times the pattern
$$
0, x_1, \ldots, x_{k-1}
$$
can be translated by an arbitrary element of $G$ to be contained in $E$.
In particular $N_{E,2}$ determines the difference multiset $E-E$ of $E$.
The $k$-deck may also be defined on an arbitrary 
locally compact abelian group,
provided we replace the summation in the definition above with integration
with respect to Haar measure.

As our primary interest is in indicator functions, we will
mainly 
consider nonnegative functions $f$. 
Another reason for considering only real functions is to avoid the extra
complication due to the fact that the functions $f$ and $\omega f$ have the
same $k$-deck whenever $\omega$ is a $k$-th root of unity.

It is evident that the functions $f(x)$ and $f_t(x) = f(x-t)$ have the same
$k$-decks for all values of $k$.
The problem we discuss in this paper is the following:
\begin{quotation}
Is the function $f:G \to \RR^{+}$ determined up to translation if we know its $k$-deck?
What if the same question is asked for indicator functions?
\end{quotation}
It is not hard to see that for $k=2$, and even for indicator functions, the
answer is negative. Indeed, suppose that we have two sets $A, B \subseteq G$
such that $-B$ is not a translate of $B$ and suppose also that the multisets
$E=A+B$ and $F=A-B$ are actually sets.
Take for example $G=\ZZ_{101}$ (the cyclic group of $101$ elements), $A=\Set{0,10,20,30}$
and $B = \Set{0,1,3}$.
Then it is easy to see that the sets $E$ and $F$ have the same
$2$-deck but are not necessarily translates of each other, e.g.\ in the example we mentioned.

In this paper we will restrict ourselves to finite cyclic groups and the emphasis will be on the $3$-deck or
triple correlation.
This problem is of significance in several fields of applied science, for example crystallography
and signal processing \cite{Pet}. See also \cite{JK} and the references
therein.

Our problem is most naturally studied with the use of the Fourier Transform
on $G$, defined for any function $f:G\to\CC$ as a function $\ft{f}$ on $\Gamma$, the group
of characters of $G$ (group homomorphisms into the multiplicative group $\Set{z\in\CC:\ \Abs{z}=1}$), given by 
$$
\ft{f}(\gamma) = \sum_{x\in G} f(x) \overline{\gamma(x)}.
$$
In the particular case of interest to us when $G = \ZZ_n$ is the cyclic
group of $n$ elements then its dual group $\Gamma$
is also isomorphic to $\ZZ_n$ and the FT of $f:\ZZ_n\to\CC$ is a function $\ft{f}:\ZZ_n\to\CC$ given by
$$
\ft{f}(k) = \sum_{j=0}^{n-1} f(j) \zeta_n^{-jk},\ \ \ k=0,\ldots,n-1,
$$
where $\zeta_n = \exp(2\pi i / n)$ is a primitive $n$-th root of unity.

It is easy to see that the Fourier Transform of the function $N_{f,k}:G^{k-1}\to\RR^{+}$,
the function $\ft{N_{f,k}}:\Gamma^{k-1}\to\CC$, is given by
\begin{eqnarray}\label{k-deck-ft}
\ft{N_{f,k}}(\xi_1,\ldots,\xi_{k-1})
 &=& \ft{f}(\xi_1)\cdots\ft{f}(\xi_{k-1}) \overline{\ft{f}(\xi_1+\cdots+\xi_{k-1})}\nonumber\\
 &=& \ft{f}(\xi_1)\cdots\ft{f}(\xi_{k-1})\ft{f}(-(\xi_1+\cdots+\xi_{k-1}))\ \ \mbox{since $f$ is real}.
\end{eqnarray}
This implies that 
\begin{equation}
\label{translation}
N_{f,k} \equiv N_{g,k} \ \Longleftrightarrow\ 
\left(\xi_1+\ldots +\xi_k=0 \Longrightarrow 
\ft{f}(\xi_1)\cdots\ft{f}(\xi_{k})=\ft{g}(\xi_1)\cdots\ft{g}(\xi_{k})\right).
\end{equation}
In particular, if $N_{f,k} \equiv N_{g,k}$ for two 
nonnegative
functions $f$ and $g$ on $G$, we immediately get $\ft{f}(0) = \ft{g}(0)$ by setting all $\xi_j = 0$.
It is also clear that $N_{f,k} \equiv N_{g,k}$ for nonnegative $f$ and $g$ implies
$N_{f,r} \equiv N_{g,r}$ for all $2\le r \le k-1$ as well, so that
identity of the $k$-decks implies the identity of all lower order $r$-decks.
Choosing $\xi_1 = -\xi_2 = \xi$ and $\xi_j = 0$ for $j\ge 3$ we get
$\Abs{\ft{f}(\xi)} = \Abs{\ft{g}(\xi)}$.
Note that if $k$ is odd then we get $\Abs{\ft{f}} \equiv \Abs{\ft{g}}$
even for two arbitrary real functions
$f$ and $g$ on $G$ with 
$N_{f,k} = N_{g,k}$ if $\ftf(0)\neq 0$ or $\ftg(0)\neq 0$.
Furthermore, if $k = 3$ and
if we know that $\ft{f}$ has no zeros on $\Gamma$ it follows
using \eqref{k-deck-ft} that the ratio $\ft{f} \bigl / \ft{g}$ is a map from $G$ to the unit circle which
is a group homomorphism, and this is equivalent to the function $f$ being a translate of the function $g$.

This reveals the fact that the main difficulty in the study of this problem is the existence of zeros in 
the Fourier Transform of the function whose $k$-deck we know.
Consider for example the case of the group $G = \ZZ_p$, $p$ a prime.
It is well known that the linear rank over $\QQ$ of the set of $p$-th roots of unity is $p-1$, and this
implies that any non-trivial $\QQ$-linear combination of at most $p-1$ such roots cannot vanish.
In other words, if we have a non-constant function $f:\ZZ_p \to \QQ$ (e.g.\ the indicator function
of a non-trivial subset of $\ZZ_p$), then its FT never vanishes on $\ZZ_p$ (which is the dual group of itself). 
By the previous discussion then the $3$-deck of any function $f:\ZZ_p \to \QQ$ determines $f$ up to translation
\cite{RS1}.

The question of whether $N_{f,k}$ determines $f$ up to translation
depends both on the group $G$ on which $f$ is defined as well as on assumed properties of $f$.
The main cases of interest are when (a) $f$ is any nonnegative function, (b) $f$ is a rational-valued function,
possibly restricted to be nonnegative, and
(c) $f$ is an indicator function. It is not hard to see, for instance, that on the group $\RR$ 
there are, for every $k$, nonnegative functions
which are not determined up to translation from their $k$-deck \cite{JK}. The same question is open if one
demands that $f$ is an indicator function of set of finite measure
although the answer is known to be positive in certain special cases of sets \cite{JK}.
On the other hand even the $3$-deck determines a function $f \in L^1(\RR)$ if it is of compact support \cite{JK}.

\subsection{Previous results}
\label{sec:previous}

In the case of cyclic groups the most significant work is that of Gr\"unbaum and Moore \cite{GM}.
This work seems largely to have gone unnoticed in the mathematical literature
although it solves the most important cases of the problem
for cyclic groups.
This is probably due to the fact that it was published in a Crystallography journal.
The following is a summary of the results in \cite{GM} regarding reconstructing $f$ on $\ZZ_n$
from its $k$-deck.
Notice that in \cite{GM} it is assumed at the outset that all functions to be reconstructed
from their $k$-deck have a non-zero sum over the group.
\begin{enumerate}
\item
For any $n$, if $f$ and $g$ are rational-valued functions on $\ZZ_n$ with the same $6$-deck then they are translates of
each other \cite[Theorem 4]{GM}.
\item
If $n$ is even 
and at least $30$ then there are sets $E, F \subseteq \ZZ_n$ which have the same $3$-deck
but are not translates of each other \cite[\S 5.3]{GM}.
\item
If $n$ is odd, $f$ and $g$ are rational-valued functions on $\ZZ_n$ with the same $3$-deck
and $\ft{f}(1) \neq 0$ then $f$ and $g$ are translates of each other \cite[Theorem 3]{GM}.
This is heavily based on a result of Lenstra \cite{L} (see our \S\ref{sec:positive}).
\item
For any $n$ suppose that
$f$ is a rational-valued function on $\ZZ_n$ and $E \subseteq \ZZ_n$, $g=\chi_E$.
Then if $f$ and $g$ have the same $4$-deck and $\ft{f}(1) \neq 0$ it follows that
$f$ and $g$ are translates of each other \cite[Theorem 5]{GM}.
It is suggested in \cite{GM} that the condition $\ft{f}(1) \neq 0$ may be unnecessary.
\item
There is no value of $k$ such that for all $n$ the equality of the $k$-deck
of two {\em real} functions $f$ and $g$ on $\ZZ_n$ implies that they are
translates of each other \cite[\S 8.2]{GM}.
\item
If $n=pqr$, with $p$ and $q$ distinct primes, and $r>1$ is an integer then there
are two rational-valued functions $f$ and $g$ on $\ZZ_n$ which have the same $3$-deck,
satisfy $\ft{f}(1) = \ft{g}(1) = 0$, and are not translates of each other
\cite[\S 5.2]{GM}.
\end{enumerate}

Radcliffe and Scott \cite{RS2} study the problem for infinite subsets of $\RR$ which are subject
to some sort of ``local finiteness'' and prove reconstructibility from the $3$-deck.
In \cite{RS1} the same authors prove reconstructibility up to translation from the $3$-deck in $\ZZ_p$,
$p$ a prime, show that almost all subsets of $\ZZ_n$ are determined up to translation by their $3$-deck and
show that any set in $\ZZ_n$ 
is determined up to translation by its $k$-deck with
$k$ being 9 times the number of distinct prime factors of $n$.

Pebody, Radcliffe and Scott \cite{PRS} study a variation of the problem. They prove that any finite
subset $E$ of the plane can be reconstructed up to rigid motion if one knows for any subset $A$ of the plane
of up to 18 points how many rigid-motion copies of $A$ are to be found in $E$.

Jaming and Kolountzakis \cite{JK} study the problem both in the case of the group $\RR$ and in
cyclic groups.
In the case of $\RR$ it is pointed out that several conditions which guarantee some sort of
analyticity of $\ft{f}$ are enough to imply that the $3$-deck of $f$ determines $f$ up to translation.
It is shown that for every $k$ there exist two nonnegative, smooth $f, g \in L^1(\RR)$ with the same
$k$-deck, which are not translates of each other. In fact for some such $f$ there exist even uncountably many,
translation inequivalent, functions $g$ which have the same $k$-deck as $f$.

It is also proved in \cite{JK} that if $E\subseteq\RR$ has finite measure, $g \in L^1(\RR)$ is nonnegative
and $\chi_E$ and $g$ have the same $3$-deck, then $g$ is itself an indicator function.
Although it is still an open problem whether any $E \subseteq \RR$ of finite measure is determined
up to translation from its $3$-deck, it is proved in \cite{JK} that if $E$ is an open set with gaps
bounded below (write $E$ as a disjoint collection of open intervals and look at the gaps so defined) then
$E$ is determined from its $3$-deck.

In the case of the cyclic group $\ZZ_n$ it is proved in \cite[Theorem 3.1 and following Remark]{JK} that
when $n=p^\alpha$, $p$ a prime larger than $2$, then the $3$-deck of a set in $\ZZ_n$ determines the
set up to translation. It is also shown that if $n=2^\alpha$ then the $4$-deck of a set determines the set
up to translation (and this is mistakenly attributed to \cite{GM}).
In \cite[Theorem 3.2]{JK} it is erroneously claimed that if $n=pq$ with $p$ and $q$ two distinct primes
then the $3$-deck is enough to reconstruct a set in $\ZZ_n$.
Given the results of \cite{GM} summarized above,
the condition $p,q>2$ clearly needs to be added and then the theorem is correct.
A corrected proof is given in our \S\ref{sec:cyclic}.
The attempt, at the end of \cite{JK}, to explain the examples given by Gr\"unbaum and Moore for the case
$n=pqr$ (see summary above) is also erroneous.

Pebody \cite{P} 
defines $r(G)$ (resp.\ $r_{\mbox{set}}(G)$) the minimum $k$ such that
the $k$-deck of a 
{\large \rm nonnegative}
rational-valued function on $G$ (resp.\ subset of $G$)
determines the function (resp.\ set) up to translation.
Improving results for the cyclic group
of Alon, Caro, Krasikov and Roditty \cite{ACKR} and Radcliffe and Scott \cite{RS1},
Pebody, computes the number $r(G)$ for all finite abelian $G$ and his result implies $r(\ZZ_n) \le 6$.
For the cyclic groups the result had already been proved in \cite{GM}.
In particular, Pebody gets that the $3$-deck determines all nonnegative
rational valued functions up to translation on the cyclic group $\Z_n$
($n\ge 3$) if and only if $n$ is a power of an odd prime or the product
of two odd primes.

\subsection{New results}
\label{sec:new}

In \S\ref{sec:cyclic} we complete the characterization of those finite cyclic
groups in which the $3$-deck determines any subset up to translation.
We show that
\begin{enumerate}
\item
If $n=p^2q$, with $p$ and $q$ distinct odd primes then any subset of $\ZZ_n$ can
be determined up to translation from its $3$-deck (Theorem \ref{posodd}).
\item
The same is true if $n=pqr$,
with $p,q,r$ distinct odd primes (Theorem \ref{posodd}).
\item
If $n=pqrd$, with $p,q$ distinct primes and $r,d>1$, then there are
two subsets $E$ and $F$ of $\ZZ_n$ with the same $3$-deck which are
not translates of each other (Theorem \ref{pqrd}).
\item
If $n=2k$, $k\ge 6$, we give two subsets $E$ and $F$ of $\ZZ_n$, not translates
of each other, which have the same $3$-deck (Theorem \ref{2k}).
This result subsumes the above mentioned result of \cite{GM} (for even $n$, $n\ge 30$)
and, we think, our examples are much easier to understand.
 
If $n$ is even and at most 10 we show that there are no such examples 
(Proposition~\ref{smalleven}).
\end{enumerate}

Thus we get the following.

\bcor\label{summary}
Every subset of the cyclic group $\Z_n$ can be determined up to translation
from its $3$-deck if and only if $n$ is a power of an odd prime or $n$ is the 
product of at most three (not necessarily distinct) odd primes or 
$n\in\{2, 4, 6, 8, 10\}$. 
\ecor

\begin{rrr}
As we were finishing this paper we came across a manuscript
by Pebody \cite{P2} where the cases of odd $n$ for which the $3$-deck
is sufficient are also determined.
Our work was done independently and, apparently, almost simultaneously.
\end{rrr}

Comparing Corollary~\ref{summary}
to the last mentioned special case of the result
of Pebody \cite{P} in the previous subsection,
we observe that the analogous characterization of the ``good'' values of $n$ is 
different if we consider nonnegative rational valued functions instead
of subsets.

Key to our results are theorems which significantly restrict the zero set of
the Fourier Transform of indicator functions of subsets of certain cyclic groups.
See for instance Lemma~\ref{key}.

In \S\ref{sec:exceptions} we study the number of subsets of $\ZZ_n$ which are not determined by
their $3$-deck up to translation.
In \cite{RS1} Radcliffe and Scott had already shown that this number is  $O(2^n/\sqrt n)$ as $n\to\infty$.
We show that this number is in fact much smaller, namely
$O(2^{-C_\epsilon n^{1-\epsilon}} 2^n)$,
for any fixed $\epsilon>0$ (Theorem \ref{th:exceptions}).


\section{For which cyclic groups the $3$-deck determines a set up to translation}
\label{sec:cyclic}

\subsection{Positive results}\label{sec:positive}

The main result of this subsection is the following:

\begin{theorem}\label{posodd}
Let $n$ be a power of an odd prime or the product of at most three
(not necessarily distinct) odd primes.
Then every subset of $\Z_n$ is uniquely determined up to translation
by its $3$-deck.
\end{theorem}

For completeness and because it needs no extra effort, 
our proof will cover not only the new results
but also the known ones. We shall use only the same theorem of
H. W. Lenstra that was used by Gr\"unbaum and Moore in \cite{GM}:

{\bf Lenstra's Theorem} \cite{L}.
{\it
If $N$ is an odd integer,
and $m$ and $N$ are coprime then 
there exists a finite sequence 
$x_1,\ldots,x_l$ of relative primes to $N$
such that $x_1=1$, $x_l=m$ and
every member except the first is the sum or difference of two 
not necessarily different previous members of the sequence.}

As we saw in the Introduction, the $3$-deck determines a nonnegative
function up to translation if its Fourier Transform has no zero. 
We shall show that if 
$n$ is a power of an odd prime or $n$ is the product of at most three
(not necessarily distinct) odd primes 
then the support of the Fourier Transform of a characteristic function on
$\Z_n$ is always rich enough to get the same conclusion. 
Our method can be considered as a generalization of the methods in
\cite{GM} and \cite{JK}.

\bd\label{extdom}
We say that $A\su\Z_n$ is an 
\emph{extendable domain} if for every 
$h:A\to\R/\Z$ additive 
(by which we mean that $h(x+y)=h(x)+h(y)$ whenever $x,y,x+y\in A$)
function there exists an $L\in\R$ such that 
$h(k)=Lk$ (mod $1$) for every $k\in A$.
\ed

\bl \label{domain}
(1) If $f$ and $g$ are nonnegative functions on $\Z_n$ with the same $3$-deck
and $\sftf$ is an extendable domain then $f$ and $g$ are translates of
each other.

(2) If $f$ and $g$ are real valued functions on $\Z_n$ with the same $3$-deck,
$\ftf(0)\neq 0$ or $\ftg(0)\neq 0$
and $\sftf$ is an extendable domain then $f$ and $g$ are translates of
each other.
\el

\begin{proof}
Suppose that $f$ and $g$ 
satisfy the conditions of (1) or (2).
We saw in the Introduction that in these cases having the same $3$-deck implies
that $\ftf$ and $\ftg$ have the same modulus.
Hence there exists a function
$h:\sftf\to\R/\Z$ such that $\ftg(l)=e^{2\pi i h(l)} \ftf(l)$.
Substituting this to (\ref{translation}) we get that
$h$ must be additive as defined in Definition \ref{extdom}.
Then, since $\sftf$ is an extendable domain, 
$h$ must be linear, thus $e^{2\pi i h(l)}$
is the restriction of a character to $\supp\ft{f}$, and so $f$ and $g$ are translates of
each other.
\end{proof}

Therefore, to prove Theorem~\ref{posodd} it is enough to prove the
following:

\begin{proposition}\label{odddomain}
If $n$ is a power of an odd prime or $n$ is the product of at most three
(not necessarily distinct) odd primes 
then the support of the Fourier Transform of a characteristic function on
$\Z_n$ is always an extendable domain.
\end{proposition}

To prove this proposition we need several facts and lemmas, some of which may be
known and/or interesting in themselves. 
The following five facts are surely known
but for completeness, and because it is easier to prove them than to find
them in the literature, we present their proofs.

\begin{notation} \label{notation}
Let $(k,l)$ denote the greatest common divisor of $k$ and $l$.
For $a|n$ let
$$
\class{a} =\{k\in\ZZ_n : (k,n)=a\} \quad \textrm{and}
$$
$$
a\Zn=\left\{0,\ a,\ 2a,\ \ldots\ ,\ \left(\frac na -1\right)a\ \right\}\su \ZZ_n.
$$

\en

\begin{fact}\label{classes}
If $f:\ZZ_n\to\Q$ then $\supp \ft{f}$ is the union of sets of the form $\class{a}$
($a|n$).
\end{fact}
\begin{proof}
We can write $\ft{f}(k) = \sum_{j=0}^{n-1} f(j) \zeta_{n,k}^j$, 
where
$\zeta_{n,k} = e^{-2\pi i k/n}$ is the $k$-th root of unity of order $n$.
The right hand side is a rational polynomial evaluated at the roots of unity.
It is well known that $\zeta_{n,k}$ is an algebraic conjugate of $\zeta_{n,l}$ over $\QQ$ if and only if $(n,k) = (n,l)$.
\end{proof}

\begin{fact}\label{periodic}
For $a|n$
a function $f:\ZZ_n\to\C$ is $a$-periodic if and only if 
$\supp \ft{f}\su\frac{n}{a}\Zn$.
\end{fact}
\begin{proof}
The space of $a$-periodic functions on $\ZZ_n$ has dimension $a$ and it is clearly spanned by the characters $\chi_l(j) = e^{2\pi i \frac{ln}{a}j}$, $l=0,1,\ldots,a-1$.
\end{proof}

\begin{fact}\label{persum}
For a function $f:\ZZ_n\to\C$ we have 
$\supp \ft{f}\su a_1\Zn\cup\ldots\cup a_k\Zn$ if and only if 
$f$ can be written in the form $f=f_1+\ldots+f_k$, 
where $f_j:\ZZ_n\to\C$ is $\frac{n}{a_j}$-periodic ($j=1,\ldots,k$).
\end{fact}
\begin{proof}
The splitting $f=\sum f_j$ is accomplished by arbitrarily splitting
$\ft{f} = \sum \ft{f_j}$, in a way that $\ft{f_j}$ is supported on $a_j\ZZ_n$,
inverting the Fourier Transform and using Fact~\ref{periodic}.
\end{proof}

\begin{fact}\label{fact:supposedly-obvious}
If $n$ is odd then any integer $k$ can be written as $k=a+b$
where $(a,n) = (b,n) = 1$.
\end{fact}
\begin{proof}
We can clearly assume that $n$ is squarefree; 
that is, it is of the form $n=p_1\cdots p_r$, where
$p_1,\ldots,p_r$ are distinct primes.
For each $j=1,\ldots,r$ let $a_j=2$ and $b_j=-1$ if $k=1\ (\md p_j)$
and let $a_j=1$ and $b_j=k-1$ otherwise. 
By the Chinese Remainder theorem 
there exist $a$ and $b$ such that
$a=a_j\ (\md p_j)$ and $b=b_j\ (\md p_j)$ for $j=1,\ldots,r$. 
Now $a+b=k\ (\md p_j)$ for each $j=1,\ldots,r$, so $a+b=k\ (\md n)$.
By choosing $a$ properly, by which we mean that we add a multiple of $n$
to $a$ if necessary, we can guarantee that $a+b=k$. 
Since each $p_j>2$, we have $a_j, b_j\neq 0\ (\md p_j)$, so $(a,n)=(b,n)=1$.
\end{proof}

\begin{fact}\label{fact:relations}
If there are two equivalence relations on a set such that both contain
at least two classes then there exist two elements which are inequivalent
w.r.t.\ both relations.
\end{fact}
\begin{proof}
If not then any two elements which are inequivalent w.r.t.\ the first relation should
be equivalent w.r.t.\ the second. This easily implies that there is only one equivalence class
w.r.t.\ the second relation, a contradiction.
\end{proof}

\begin{lemma}\label{aa}
If $a$ is a divisor of the odd $n$ and $\class{a}\su A\su a\Zn$
then $A$ is an extendable domain.
\end{lemma}

\begin{proof}
Let $h:A\to\R/\Z$ be an additive function. 
It is enough to prove that 
\begin{equation}\label{claima}
m\in\Z, ma\in A \Longrightarrow h(ma)=m h(a)\ (\md 1),
\end{equation}
since then for any $L\in\R$ such that $h(a)=La\ (\md 1)$
we get that $h(ma)=mh(a)=Lma\ (\md 1)$, 
which is exactly what we want to show.

Let $N=n/a$. Note that $ma\in\class{a}$ holds if and only if $m$ and $N$
are coprime.
Using the previously stated Lenstra's Theorem 
for $N=n/a$ and an $m$ such that $m$ and $N$ are coprime
we get a sequence 
$x_1,\ldots,x_l$ of relative primes to $N$
such that $x_1=1$, $x_l=m$ and
every member except the first is the sum or difference of two 
not necessarily different previous members of the sequence.
Note that, since $x_i$ and $n/a$ are coprime, $x_i a\in\class{a}\su A$
for each $i$. 
Then by induction we get that $h(x_i a)=x_i h(a)$, and so 
(\ref{claima}) holds whenever $m$ and $n/a$ are coprime.
Then (\ref{claima}) in
the general case follows by using Fact \ref{fact:supposedly-obvious}.
\end{proof}

\begin{corollary}\label{cor}
If $n$ is odd, $f:\ZZ_n\to\Q$ and $a\in\sftf\su a\Zn$ then
$\sftf$ is an extendable domain.

In particular, the following two statements hold:
\begin{itemize}
\item[(i)] If $n$ is odd, $f:\ZZ_n\to\Q$ and $\ftf(1)\neq 0$
then $\sftf$ is an extendable domain.
\item[(ii)]
If $n$ is a power of an odd prime and  $f:\ZZ_n\to\Q$
then $\sftf$ is an extendable domain.
\end{itemize}
\end{corollary}

\bp
The first statement follows immediately from  Lemma~\ref{aa} and
Fact~\ref{classes}.
If $a=1$ then we get (i). Statement (ii) is also a special case of
the first statement since if $n=p^k$ and $l$ is minimal such that
$p^l\in\sftf$ then $p^l\in\sftf\su p^l\Z_n$.
\ep

\begin{lemma}\label{0}
If $\chi_E=f_a+f_b$, and $f_a$ and $f_b$  are periodic $\Z\to\C$
functions with coprime periods $a$ and $b$ then $\chi_E$ is periodic
with period $a$ or $b$. 
\el

\begin{proof} 
Using the periodicity and $\chi_E=f_a+f_b$, for any $k,n,l\in\N$ we get
\begin{equation}\label{range}
f_b(ak+bn+l)=f_b(ak+l)=\chi_E(ak+l)-f_a(ak+l)=\chi_E(ak+l)-f_a(l).
\end{equation}
Since $a$ and $b$ are coprime, for any fixed $l\in\Z$ every integer
can be written in the form $ak+bn+l$. 
Thus (\ref{range}) implies that for any fixed $l\in\Z$ the range of 
$f_b$ is a subset of $\{-f_a(l),1-f_a(l)\}$. 
This implies that $f_a$ or $f_b$ must be constant on $\Z$.
\end{proof}

\begin{lemma}\label{na}
Suppose that $a$ and $b$ are divisors of $n$, $n/a$ and
$n/b$ are coprime, $E\su\ZZ_n$, and
$\sFTE\su a\Zn \cup b\Zn$.
Then $\sFTE \su a\Zn$ or $\sFTE \su b\Zn$.
\el

\begin{proof}
By Fact~\ref{persum}, $\sFTE\su a\Zn \cup b\Zn$ implies the existence
of an $n/a$-periodic function $f$ and an $n/b$ periodic function $g$
such that $\chi_E=f+g$. Since $n/a$ and
$n/b$ are coprime, by Lemma~\ref{0}, we get that
$\chi_E$ must be $n/a$-periodic or $n/b$ periodic.
By Fact~\ref{periodic}, this implies that 
$\sFTE \su a\Zn$ or $\sFTE \su b\Zn$.
\end{proof}

\begin{lemma}\label{ab}
Suppose that $a$ and $b$ are coprime divisors of $n$, $E\su\ZZ_n$, 
$$
\sFTE\su (a\Zn \cup b\Zn \sm ab\Zn) \cup \Set{0}. 
$$
Then $\sFTE \su a\Zn$ or $\sFTE \su b\Zn$.
\el

\begin{proof}
Let $c=\frac{n}{ab}$. 
By Fact~\ref{persum}, $\sFTE\su a\Zn \cup b\Zn$ implies that $\chi_E$
can be written in the form 
$\chi_E=f_a+f_b$, where $f_a$ is $bc$-periodic and $f_b$ is $ac$-periodic.
Applying Lemma~\ref{0} we get that for each $t=0,1,\ldots,p-1$ the function
$\chi_E(kc+t)$ is $a$-periodic or $b$-periodic.
Let $m_t$ be the number of points of the form $kc+t$ in $E$.
Then
\begin{multline}\label{mt}
\textrm{$m_t\in\{0,1,\ldots,ab\}$ is divisible by $a$ or $b$;}\\ 
\textrm{by $a$ if $\chi_E(kc+t)$ is $b$-periodic and 
by $b$ if $\chi_E(kc+t)$ is $a$-periodic.}
\end{multline}

A straightforward calculation shows that for any $s\in\Z$
$$
\FTE(sab)=\sum_{t=0}^{c-1} m_t \left(e^{\frac{2\pi i s}{c}}\right)^t.
$$
Since we assumed that $\FTE(sab)=0$ for $s=1,\ldots,c-1$, 
we get that the $c-1$ $c$-th roots of unity $e^{\frac{2\pi i s}c}$
($s=1,\ldots,c-1$) are all roots of the $(c-1)$-th order polynomial
$\sum_{t=0}^{c-1} m_t z^t$.
Hence $\sum_{t=0}^{c-1} m_t z^t$ must be a constant multiple of
$\Pi_{s=1}^{c-1}(z-e^{\frac{2\pi i s}c})=\sum_{t=0}^{c-1} z^t$ and so
all $m_t$ must be the same. 
Using (\ref{mt}) and that $a$ and $b$ are coprime this implies that 
$\chi_E(kc+t)$ is $a$-periodic for each $t$ or $b$-periodic for each $t$.
Thus $\chi_E$ is $ac$-periodic or $bc$-periodic hence, by Fact~\ref{periodic},
$\sFTE \su a\Zn$ or $\sFTE \su b\Zn$.
\end{proof}

\bl \label{linear}
If $a$ and $b$ are coprime divisors of the odd $n$ and
$\class{a} \cup \class{b} \cup \{ab\} \su A \su a\Zn \cup b\Zn$
then $A$ is an extendable domain. 
\el 

\begin{proof}
Let $h:A\to\R/\Z$ be an additive function. 
We have to find an $L\in\R$ such that 
$h(k)=Lk$ (mod $1$) for every $k\in A$.

By Lemma~\ref{aa}, $A\cap a\Zn$ and $A\cap b\Zn$ are extendable domains,
so there exist $L_a$ and $L_b$ such that 
\begin{equation}\label{Lab}
h(k)=L_a k\ (\md 1) \quad \textrm{if } k\in A\cap a\Zn,
\qquad \textrm{and} \qquad
h(k)=L_b k\ (\md 1) \quad \textrm{if } k\in A\cap b\Zn.
\end{equation}
Note that for any $u,v\in\Z$, $L_a$ can be replaced by $L_a+\frac ua$
and $L_b$ can be replaced by $L_b+\frac vb$ in (\ref{Lab}).
Thus it is enough to find $u,v\in\Z$ such that 
$L_a+\frac ua=L_b+\frac vb$, which is equivalent to
\begin{equation}\label{uv}
ub-va=L_a ab - L_b ab.
\end{equation}

Using $ab\in A$ and (\ref{Lab}), we get $L_a ab=h(ab)=L_b ab\ (\md 1)$,
so $L_a ab - L_b ab\in\Z$. 
Then, since $a$ and $b$ are coprime,
there exists $u,v\in\Z$ for which (\ref{uv}) holds, 
which completes the proof.
\end{proof}

In the sequel we shall use Fact~\ref{classes} in the proofs 
many times without explicitly citing it.

\bl\label{pqd}
Let the odd $n=pqd$, where $p$ and $q$ are two distinct primes, and $d$ is a prime
or $d=1$. 
If $E\su\ZZ_n$ and $\sFTE\su p\Zn\cup q\Zn$ then 
$\sFTE$ is an extendable domain.
\el

\begin{proof}
If $p,q\in\sFTE$ then, applying Lemma~\ref{ab} for $a=p$, $b=q$, we
get that
$\sFTE\cap pq\Zn\neq\Set{0}$,
which implies that $pq\in\sFTE$.
Then we can apply Lemma~\ref{linear} to get that $\sFTE$ is indeed an
extendable domain.

So we can suppose by symmetry that $q\not\in\sFTE$. 
Then $\sFTE\su p\Zn\cup q\Zn$ implies that 
$\sFTE\su p\Zn\cup qd\Zn$.
Then, in case of $d\neq p$ by Lemma~\ref{na}, in case of $d=p$ clearly,
we have 
$\sFTE\su p\Zn$ or $\sFTE\su qd\Zn$. 

If $\sFTE\su qd\Zn$ then $\sFTE=\class{qd}\cup\{0\}$ or
$\sFTE=\{0\}$, so we are done by Lemma~\ref{aa}.

So we can suppose that $\Set{0}\neq\sFTE\su p\Z_n$.
If $p\in\sFTE$ then by Corollary~\ref{cor} $\sFTE$ is an extendable domain,
so we can suppose that $p\not\in\sFTE$. 
Then, by Fact~\ref{classes}, $\sFTE$ can be only  
$\class{pd}\cup\Set{0}$ or $\class{pq}\cup\Set{0}$ or
$\class{pq}\cup\class{pd}\cup\Set{0}$ with $d\neq 1,q$.
The last case is impossible by Lemma~\ref{na}
(for $a=pq, b=pd$), 
while in the first two cases Lemma~\ref{aa} implies
that $\sFTE$ is an extendable domain.
\end{proof}

The following lemma about the possible support of the 
Fourier Transform of characteristic functions on $\Zn$
is the key for handling the hardest case when $n$ is
the product of three distinct primes. 
This statement might be useful in other applications, too.
\bl \label{key}
Suppose $p, q$ and $r$ are pairwise coprime, but not necessarily primes.
Let $n=pqr$ and let
$E\su\ZZ_n$. Then
$$
p,q\in\sFTE\su p\Zn \cup q\Zn \cup r\Zn
\Longrightarrow
\left(\exists z\in\{1,2,\ldots,r-1\}\right)\  zpq\in\sFTE.
$$
\el

\begin{proof}
Suppose that 
for each $z=1,2,\ldots,r-1$ we have
$$
0=\FTE(zpq)=\sum_{c=0}^{r-1}\sum_{k=0}^{pq-1}
\chi_E(kr+c)\e{(kr+c)zpq}{pqr}=
\sum_{c=0}^{r-1}\left(\sum_{k=0}^{pq-1} \chi_E(kr+c) \right)
\left(\e{z}{r}\right)^c.
$$

This implies that $\sum_{k=0}^{pq-1} \chi_E(kr+c)$ must be the same
for each $c\in\ZZ_r$; 
that is, 
\begin{equation}\label{constant}
\sum_{k=0}^{pq-1} \chi_E(kr+c_1)-\chi_E(kr+c_2)=0
\qquad (c_1, c_2\in\Z_r).
\end{equation}

For each $j\in\ZZ_n$ ($n=pqr$) let 
$(a_j,b_j,c_j)\in \Z_p \times \Z_q \times \Z_r$ 
be the unique triple such that 
$$
j=a_j\ \md p, \qquad j=b_j \ \md q \quad \tand \quad j=c_j \ \md r,
$$
and let $\phi$ be the inverse of the above 
$\ZZ_n\to\Z_p \times \Z_q \times \Z_r$ bijections; 
that is, $\phi(a,b,c)$ ($a\in\Z_p, b\in\Z_q, c\in\Z_r$)
is the unique element of $\ZZ_n$ for which 
$$
\phi(a,b,c)=a \ \md p, \qquad \phi(a,b,c)=b \ \md q 
\quad \tand \quad \phi(a,b,c)=c \ \md r.
$$

Since $\sFTE\su p\Zn \cup q\Zn \cup r\Zn$, $\chi_E$ can be written
as $\chi_E=f+g+h$, where $f$ is $qr$-periodic, $g$ is $pr$ periodic
and $h$ is $pq$-periodic. 

Since $f$ is $qr$-periodic, $f(\phi(a,b,c))$ does not depend on $a$,
so  $f\circ \phi$ can be written in the form $f(\phi(a,b,c))=F(b,c)$. 
Similarly $g\circ \phi$ and $h\circ \phi$ can be written as
$g(\phi(a,b,c))=G(a,c)$ and $h(\phi(a,b,c))=H(a,b)$. 
So using the notation $E'=\phi^{-1}(E)\su\Z_p \times \Z_q \times \Z_r$
we get that
$$
\chi_{E'}(a,b,c)=F(b,c)+G(a,c)+H(a,b) \qquad a\in\Z_p, b\in\Z_q, c\in\Z_r.
$$

We claim that there exist $c_1,c_2\in\Z_r$ such that 
neither the $\Z_q\to\R$ function $F(\cdot\,,c_1)-F(\cdot\,,c_2)$, nor the 
$\Z_p\to\R$ function $G(\cdot\,,c_1)-G(\cdot\,,c_2)$ is constant.
Indeed,
$F(\cdot\,,c_1)-F(\cdot\,,c_2)$ being constant defines an equivalence relation on $\ZZ_r$
and so does $G(\cdot\,,c_1)-G(\cdot\,,c_2)$ being constant.
If there is only one equivalence class w.r.t.\ the first relation then
$F$ can be written as $F(b,c)=u(b)+v(c)$ which implies
$$
\chi_E(j)=\chi_{E'}(a_j,b_j,c_j)=v(c_j)+G(a_j,c_j)+u(b_j)+H(a_j,b_j),
$$
and this would in turn imply $\FTE(p)=0$, contradicting our assumption.
Hence there are at least two classes w.r.t.\ the first relation.
Similarly there are two classes w.r.t.\ the second relation and
using Fact \ref{fact:relations} we obtain our claim.

On the other hand, since
$$
(F(b,c_1)-F(b,c_2))+(G(a,c_1)-G(a,c_2))=
\chi_{E'}(a,b,c_1)-\chi_{E'}(a,b,c_2)\in\{-1,0,1\}
$$
for any $a\in\Z_p, b\in\Z_q$,
we have
$$
\Range(F(\cdot\,,c_1)-F(\cdot\,,c_2)) + 
\Range(G(\cdot\,,c_1)-G(\cdot\,,c_2)) \su \{-1,0,1\}.
$$
Since by the previous paragraph 
$\Range(F(\cdot\,,c_1)-F(\cdot\,c_2))$ 
and $\Range(G(\cdot\,,c_1)-G(\cdot\,,c_2))$ have
at least two elements, this implies that they must be of the
form 
$$
\Range(F(\cdot\,,c_1)-F(\cdot\,,c_2))=\{A,A+1\},
$$ 
$$ 
\Range(G(\cdot\,,c_1)-G(\cdot\,,c_2))=\{-A,-A-1\}
$$
for some $A\in\R$.

Let $l_1\in\{1,\ldots,q-1\}$ be the number of elements $b\in\Z_q$ 
for which $F(b,c_1)-F(b,c_2))=A$ and 
$l_2\in\{1,\ldots,p-1\}$ be the number of elements $a\in\Z_p$ for which 
$G(a,c_1)-G(a,c_2))=-A$.

Then, combining this with (\ref{constant}) we get
\begin{eqnarray}
0 &=& \sum_{k=0}^{pq-1} \chi_E(kr+c_1)-\chi_E(kr+c_2)\nonumber\\
&=& \sum_{a\in\Z_p}\sum_{b\in\Z_q} \chi_{E'}(a,b,c_1)-\chi_{E'}(a,b,c_2)\nonumber\\
&=& \sum_{a\in\Z_p}\sum_{b\in\Z_q} F(b,c_1)-F(b,c_2)+G(a,c_1)-G(a,c_2)\nonumber\\
&=& pl_1 A + p(q-l_1)(A+1) + ql_2(-A)+q(p-l_2)(-A-1)\nonumber\\
&=& -l_1p+l_2q\nonumber,
\end{eqnarray}
which is a contradiction since $l_1 p$ cannot be divisible by $q$.
\end{proof}

\bl \label{pqrsupp}
Let $n=pqr$ with $p,q,r$ three distinct primes, 
$E\su\ZZ_n$. Then
$$
p,q,r\in\sFTE\su p\Zn \cup q\Zn \cup r\Zn
\Longrightarrow
\sFTE = p\Zn \cup q\Zn \cup r\Zn.
$$
\el

\begin{proof}
Suppose that 
$p,q,r\in\sFTE\su p\Zn \cup q\Zn \cup r\Zn$.
By Lemma~\ref{key} we have 
$\class{pq}\cup\class{pr}\cup\class{qr}\su\sFTE$.
Since $E$ cannot be empty, $\class{pqr}\su\sFTE$. 
Since $p\Zn \cup q\Zn \cup q\Zn=\class{p}\cup\class{q}\cup\class{r}
\cup\class{pq}\cup\class{pr}\cup\class{qr}\cup\class{pqr}$, 
this completes the proof.
\end{proof}

\bl \label{pqrdomain}
If $a,b,c\in\Z$ are pairwise coprime then $a\Zn\cup b\Zn \cup c\Zn$ is
an extendable domain.
\el

\begin{proof}
Let $A=a\Zn\cup b\Zn \cup c\Zn$ and let $h:A\to\R/\Z$ be additive.
Then it is easy to show that for $\al=h(a)/a, \be=h(b)/b, \ga=h(c)/c$ we have
\begin{equation}\label{m}
h(ma)=\al ma, \quad h(mb)=\be mb, \quad h(mc)=\ga mc \quad (\md 1)
\quad (m\in\Z).
\end{equation}

It is enough to find $u,v,w\in\Z$ such that
\begin{equation}\label{uvw}
\al+\frac ua=\be + \frac vb = \ga+\frac wc
\end{equation}
since then $h(x)=Lx\ (\md 1)$ would follow for 
$L=\al+\frac ua=\be + \frac vb = \ga+\frac wc$  from (\ref{m}).

For $v\in\Z$ there exist $u$ and $w$ such that (\ref{uvw}) holds if
$$
va=\al ab-\be ab\ (\md b) \quad \tand vc=\ga bc -\be bc\ (\md b),
$$
which hold for some $v\in\Z$ if and only if
\begin{equation}\label{enough}
\al ab-\be ab\in\Z, \quad \ga bc -\be bc\in\Z \ \tand \ 
c(\al ab-\be ab)=a(\ga bc -\be bc)\ (\md b).
\end{equation}

Using (\ref{m}) for $m=a,b,c$ we get that
$$
\al ab=\be ab, \qquad \be bc=\ga bc 
\quad \tand \quad \al ac=\be ac  \qquad (\md 1),
$$
which implies (\ref{enough}) and so completes the proof.
\end{proof}

\begin{proof}(Proposition~\ref{odddomain})
By Corollary \ref{cor} we are done if $n$ is a power of an odd prime
or if $1\in\sFTE$. 
So we can suppose that $1\not\in\sFTE$ and $n=pqr$,
where $p$ and $q$ are different primes and $r$ is a prime or $r=1$.
If $r=1$ or $r=p$ or $r=q$ then $n$ equals $pq$ or $p^2q$ or $pq^2$
and so $1\not\in\sFTE$ implies that $\sFTE\su p\Zn\cup q\Zn$,
hence we are done by Lemma~\ref{pqd}.

Therefore we can suppose that $1\not\in\sFTE$ and $n=pqr$, where
$p, q$ and $r$ are distinct primes. Then we have
$\sFTE\su p\Zn\cup q\Zn\cup r\Zn$.

If $p,q,r\in\sFTE$ then, by Lemma~\ref{pqrsupp}, we have
$\sFTE = p\Zn\cup q\Zn\cup r\Zn$, which is an extendable domain
by Lemma~\ref{pqrdomain}. 

Otherwise, we can suppose by symmetry that $r\not\in\sFTE$ and so
$\sFTE\su p\Zn\cup q\Zn$. Then we are done by Lemma~\ref{pqd}

This completes the proof of Proposition~\ref{odddomain} and so 
also the proof of Theorem~\ref{posodd}.
\end{proof}

If $n$ is even then we get positive results for small $n$:

\begin{proposition}\label{smalleven}
For $n=2, 4, 6, 8$ and $10$ every subset of
$\ZZ_n$ is uniquely determined up to translations by its $3$-deck.
\end{proposition}

For $n=2, 4$ and $6$ this statement follows very easily from the definition
of $3$-deck. For both $n=8$ and $n=10$ one can provide proofs using the
lemmas and the method of this section.  
However, in these cases
there are only $2^8$ and $2^{10}$ subsets of $\ZZ_n$, so one can easily check
(and we indeed did check) the statement by computer. 
Hence we omit the quite complicated detailed proof, in which many
cases have to be distinguished and no new idea is needed.

\subsection{Negative results}\label{negative}

\bt \label{2k}
Let $n=2k$ with $k\ge 6$ integer.
Then there exists $E,F\su\ZZ_n$ such that they are not translates of each other,
however they have the same $3$-deck.
\et

\begin{proof}
Let 
$$
E=\{0\}\cup\{3, 4,\ldots, k-1\} \cup \{k+1, k+2\}, \quad \textrm{and} \quad
F=\{0, 1\}\cup\{3, 4,\ldots, k-1\} \cup \{k+2\}. 
$$
Since $k\ge 6$, both $E$ and $F$ contain a unique block of $k-3$ consecutive 
numbers. Thus if a translation takes $E$ to $F$ then this block of $E$ must
be taken to the block of $F$. Since these blocks are identical, the
translation must be the identity. But $E\neq F$, so we proved that
they are not the translates of each other.

By (\ref{k-deck-ft}),
for checking that $E$ and $F$ have the same $3$-deck we have to
show that for the Fourier Transforms of their characteristic function we have 
\begin{equation}\label{toprove1}
s_1+s_2+s_3=0\ (\md 2k) \Longrightarrow
\FTE(s_1)\FTE(s_2)\FTE(s_3)=\FTF(s_1)\FTF(s_2)\FTF(s_3).
\end{equation}

Letting $z=\zeta_{2k}^{-s}=e^{-2\pi i\frac{s}{2k}}$ we have
\begin{eqnarray}\label{fte}
\FTE(s)&=&
\zz^{-0s}+\left(\zz^{-3s}+\zz^{-4s}+\ldots+\zz^{-(k-1)s}\right)
+\zz^{-(k+1)s}+\zz^{-(k+2)s}\\ \nonumber
&=&1+(z^3+z^4+\ldots+z^{k-1})+z^{k+1}+z^{k+2}=
(1-z+z^3)(1+z+\ldots+z^{k-1}),\nonumber
\end{eqnarray}
and similarly
\begin{equation}\label{ftf}
\FTF(s)=1+z+(z^3+z^4+\ldots+z^{k-1})+z^{k+2}=
(1-z^2+z^3)(1+z+\ldots+z^{k-1}).
\end{equation}
If $s$ is even but $s\neq 0\ (\md 2k)$ then 
$$
1+z+\ldots+z^{k-1}=\frac{z^k-1}{z-1}
=\frac{e^{-2\pi i\frac{ks}{2k}}-1}{e^{-2\pi i\frac{s}{2k}}-1}=0,
$$
and so $\FTE(s)=\FTF(s)=0$.

Since $s_1+s_2+s_3=0\ (\md 2k)$ implies that at least one of $s_1, s_2$ and 
$s_3$ is even, we get that (\ref{toprove1}) clearly holds unless 
at least one of $s_1, s_2$ and $s_3$ is zero. 

So suppose that at least one of $s_1, s_2$ and $s_3$ is zero.
Then, in order to check (\ref{toprove1}), we
have to show that 
\begin{equation}\label{s-s}
\FTE(s)\FTE(-s)=\FTF(s)\FTF(-s) \qquad (s\in\Z_{2k}).
\end{equation}
This is just a restatement of the fact that $E$ and $F$ have the same $2$-deck,
which is clearly true as $E$ is a translate of $-F$.
\end{proof}

\bt \label{pqrd} 
Let $n=pqrd$ with $p,q$ two distinct primes and $r,d>1$ integers. 
Then there exist $E,F\su\ZZ_n$ such that they are not translates of each other,
however they have the same $3$-deck.
\et

\begin{proof}
Let
$$
A=\left\{\frac{l_1 n}{q}+kd : 
k\in\{0,1,\ldots,r-1\},l_1\in\{0,1,\ldots,q-1\}\right\},
$$
$$
B=\left\{\frac{l_2 n}{p}+kd : 
k\in\{0,1,\ldots,r-1\},l_2\in\{0,1,\ldots,p-1\}\right\},
$$
$$
E=A\cup(B+1) \qquad \textrm{and} \qquad F=A\cup(B+d+1).
$$

Then 
$$
\FTE(s)=
\underbrace{\left(\sum_{k=0}^{r-1}\e{kds}{n}\right)}_
{0 \tif pq|s\ {\rm but}\ pqr \nmid s} 
\cdot \left(
\underbrace{\sum_{l_1=0}^{q-1}\e{l_1 s}{q}}_
{q \tif q|s, \ 0 \tif q\nmid s} 
+\e{s}{n}\cdot 
\underbrace{\sum_{l_2=0}^{p-1}\e{l_2 s}{p}}_
{p \tif p|s,\ 0 \tif p\nmid s} 
\right)
$$
and
$$
\FTF(s)=
\left(\sum_{k=0}^{r-1}\e{kds}{n}\right)\cdot
\left(\sum_{l_1=0}^{q-1}\e{l_1 s}{q} +
\e{s(d+1)}{n}\cdot \sum_{l_2=0}^{p-1}\e{l_2 s}{p}\right).
$$

Thus 
$$
\FTE(s)\neq 0 \Longleftrightarrow \FTF(s)\neq 0 \Longleftrightarrow
(p|s \tor q|s) \tand (pq\!\!\not| s \tor pqr|s),
$$
hence
$$
\sFTE=
\underbrace{\{s\in\ZZ_n: p|s, q\!\!\not| s\}}_{S_p} \cup 
\underbrace{\{s\in\ZZ_n: q|s, p\!\!\not| s\}}_{S_q} \cup 
\underbrace{\{s\in\ZZ_n: pqr|s\}}_{S_{pqr}}.
$$

For checking that $E$ and $F$ have the same $3$-deck we have to
show that 
\begin{equation}\label{toprove}
s_1+s_2+s_3=0 (\md n) \Longrightarrow
\FTE(s_1)\FTE(s_2)\FTE(s_3)=\FTF(s_1)\FTF(s_2)\FTF(s_3).
\end{equation}

We have nothing to prove unless $s_1,s_2,s_3\in\sFTE$.
So suppose that $s_1,s_2,s_3\in\sFTE$.

Note that $\FTE(s)=\FTF(s)$ unless $s\in S_p$.
Thus if none of $s_1,s_2,s_3$ are in $S_p$ then we are done.

It is impossible that exactly one of them is in $S_p$ (because of
divisibility by $q$).

If two of them, say $s_1$ and $s_2$, are in $S_p$ then $s_3$ is in $S_p$
or in $S_{pqr}$. In both cases it is easy to check (\ref{toprove}).

Finally, suppose that $F=E+t$ (mod $n$). 
Since both $E$ and $F$ have $qr$ elements
that are $0$ mod $d$ and $pr$ elements that are $1$ mod $d$, $t$ must
be of the form $t=md$.
Thus we must have $A=A+md$ and $B+d=B+md$ (mod $n$). 
But $B$ consists of blocks which are arithmetic progressions of length $r$ and step $d$,
and these are regularly spaced at intervals of length $qrd$.
Hence $B+d=B+md$ (mod $n$) can only happen if $d-md$ is a multiple of $qrd$, or, equivalently,
if $m = 1 \bmod qr$.
On the other hand, by the similar structure of $A$ it follows that $m = 0 \bmod pr$, which is
a contradiction.
\end{proof} 

\subsection{Results about real-valued functions}

Given the results we have proved so far
we can also characterize those values of $n$ for which 
the $3$-deck determines the characteristic function of any
nonempty subset of $\Zn$ up to translation 
even among all $\Zn\to\R$ functions.

\bt\label{extra}
For $n\ge 3$ the following three statements are equivalent.
\begin{itemize}

\item[(i)] $n$ is a power of an odd prime or $n$ is the product of
at most $3$ (not necessarily distinct) odd primes.

\item[(ii)] The support of the Fourier Transform of any characteristic
function on $\Z_n$ is an extendable domain.

\item[(iii)] 
If for some $\emptyset\neq E\su\Zn$ and $g:\Zn\to\R$, 
$\chi_E$ and $g$ have the same $3$-deck then they are 
translates of each other.
\end{itemize}
\et

\begin{proof}
(i)$\Rightarrow$(ii):
This is exactly Proposition~\ref{odddomain}.

(ii)$\Rightarrow$(iii):
If $E\neq\emptyset$ then $\FTE(0)\neq 0$ and so
by Lemma~\ref{domain}~(2) we get that
$\chi_E$ and $g$ are indeed
translates of each other.

(iii)$\Rightarrow$(i):
If $n$ is odd and (i) does not hold then, by Theorem~\ref{pqrd},
there exists counterexamples for (iii), even with $g$ being
a characteristic function.

Now suppose that $n>2$ is even and let $E=\{1,2,\ldots,n/2\}$.
It is easy to check that $\sFTE=\{0,1,3,5,\ldots,n-1\}$.
Let
$$
h_{\al}(l)= \begin{cases}
\al  & \text{if}\quad l=1,  \\
-\al  & \text{if}\quad l=-1,  \\
0     & \text{otherwise.}
\end{cases}
$$
Let $g_\al$ be the inverse Fourier Transform of the 
function $G_\al(l)=e^{2\pi i h_{\al}(l)}\cdot\FTE(l)$ on $\Z_n$.
Since $h_\al$ is an odd function, $G_\al(-l)=\overline{G_\al(l)}$,
and so $g_\al$ is a real valued function.
Since $h_\al$ is additive on $\sFTE$, 
the right hand side of (\ref{translation}) holds for $k=3$, $f$ and $g=g_{\al}$,
and so $N_{\chi_E,3}=N_{g_\al, 3}$.
This way we get continuum many distinct $g_\al:\Zn\to\R$ functions.
Since $\chi_E$ has only finitely many translates (iii) cannot hold for 
every $g_{\al}$.
\end{proof}

\begin{example}\upshape
Let $n\ge 4$ be arbitrary, $f=0$ on $\Zn$ and $g(k)=\cos\frac{2k\pi}n$
($k\in\Zn$).
Then clearly $\ftf=0$ and one can check that
$$
\ftg(l)=\begin{cases}
n/2  & \text{if}\quad l=1,  \\
-n/2  & \text{if}\quad l=-1,  \\
0     & \text{otherwise.}
\end{cases}
$$
Then it is easy to check that the righthand-side of
(\ref{translation}) holds for $k=3$,
so $N_{f,3}=N_{g,3}$, however $f$ and $g$ are clearly not
translates of each other.

This shows that if we allow $E=\emptyset$ in (iii) of Theorem~\ref{extra}
then (i)$\Rightarrow$(iii) is not true any more.
\end{example}

\br
It is proved in \cite{JK} (Proposition 2.7) that
if $f$ is the characteristic function of a subset of $\R$ of finite measure
and $g\in L^1(\R)$ is a nonnegative function such that $N_{f,3}=N_{g,3}$
then there $g$ must be equal to a characteristic function almost everywhere. 

One can check that the same proof works on $\Zn$ as well. This has the
following consequences.
\begin{enumerate}
\item 
The characteristic functions on $\Zn$ are determined up to
translation by their $3$-deck among nonnegative functions if and only
if they are determined up to translation among characteristic functions;
that is, by Corollary~\ref{summary}, if and only if 
$n$ is a power of an odd prime or $n$ is the 
product of at most three (not necessarily distinct) odd primes or 
$n\in\{2, 4, 6, 8, 10\}$.
\item 
Only the (at most finitely many) characteristic functions can be
nonnegative among the (continuum many) $g_{\al}$ functions
of the proof of (iii)$\Rightarrow$(i) of Theorem~\ref{extra}.
\end{enumerate}
\erem


\section{The percentage of subsets of $\ZZ_n$ not determined by
their $3$-deck up to translation}
\label{sec:exceptions}

As we mentioned in the Introduction,
in \cite{RS1} Radcliffe and Scott proved that almost all subsets of $\ZZ_n$ are determined
up to translation by their $3$-deck.
More specifically they proved that the fraction of subsets of $\ZZ_n$ whose Fourier Transform
vanishes somewhere is at most $C_\epsilon \bigl / n^{1/2-\epsilon}$, for any $\epsilon>0$,
and, since any set whose FT does not
vanish is uniquely determined from its $3$-deck, this proves that a fraction at most
$C_\epsilon \bigl / n^{1/2-\epsilon}$ of the possible sets are not determined by their $3$-deck.

Furthermore, it is easy to see that
the probability of having the FT of a random subset of $\ZZ_n$ vanish somewhere
is at least $C \bigl / \sqrt{n}$.
For this one takes $n$ to be even and examines the FT of the random set at $n/2$.
The vanishing there is equivalent to a random subset of a set of $n/2$ ones and $n/2$
minus-ones having a vanishing sum.
This probability is equal to ${n \choose n/2} \bigl / 2^n \sim C \bigl / \sqrt{n}$.

However, here we show that the probability that a random subset of $\ZZ_n$ is
not uniquely determined up to translation by its $3$-deck is exponentially small (Theorem \ref{th:exceptions}).
When talking about random sets in this section we mean that all subsets of $\ZZ_n$
are equally probable. This is the same as tossing an independent fair coin for each
element of $\ZZ_n$ to decide membership in the random set.

\begin{lemma}\label{lm:lin-indep}
Suppose $u_1, \ldots, u_m$ are vectors in a vector space $V$ and 
that the collection $u_1,\ldots,u_D$, $D\le m$, are linearly independent.
Suppose also that $\epsilon_j$, $j=1,\ldots,m$, are $\Set{0,1}$-valued random variables
which are unbiased and independent.
Then
\begin{equation}\label{prob-bound}
\Prob{\sum_{j=1}^m \epsilon_j u_j = 0} \le 2^{-D}.
\end{equation}
\end{lemma}

\begin{proof}
Since $u_1,\ldots,u_D$ are independent,
for any fixed $\epsilon_{D+1},\ldots,\epsilon_m$, the $2^D$ possible values
of $\sum_{j=1}^m \epsilon_j u_j$ are all distinct, so only at most
one of them can be zero.
\end{proof}

\begin{corollary}\label{cor:at-one}
If $E \subseteq \ZZ_n$ is random then
$$
\Prob{\ft{\chi_E}(k) = 0} \le 2^{-Cn\bigl /(k,n)\log\log n},
$$
for some absolute constant $C>0$ and for all $k\in\ZZ_n$.
\end{corollary}
\begin{proof}
Let $\omega = e^{2\pi i / n}$.
Then
\begin{equation}\label{ft-at-one}
\ft{\chi_E}(k) = \sum_{j=0}^{n-1} \epsilon_j \omega^{kj},
\end{equation}
where the $\epsilon_j$, $j=0,\ldots,n-1$, are independent, unbiased, $\Set{0,1}$-valued
random variables.

It is well known that the algebraic order of $\omega^k$ over the field $\QQ$ is $\phi(n/(k,n))$, 
where $\phi(n)$ is the Euler
function which counts how many numbers from $1$ to $n$ are coprime to $n$.
It is also well known \cite{HW} that
$\phi(n) \ge C n \bigl / \log\log n$.
This means that if $P(x)$ is a polynomial with rational coefficients
and degree $< C n\bigl /(k,n) \log\log n$ then $P(\omega^k) \neq 0$.
This, in turn, implies that the complex numbers
$$
1,\omega^k,\omega^{2k},\ldots,\omega^{(C n/(k,n)\log\log n)\cdot k}
$$
are $\QQ$-linearly independent.
Applying Lemma \ref{lm:lin-indep} to the random sum \eqref{ft-at-one}
in the vector space $\CC$ over $\QQ$ we get our result.
\end{proof}

\begin{corollary}\label{cor:odd}
If $n$ is odd then the probability that a random subset of $\ZZ_n$
is not uniquely determined by its $3$-deck is at most $2^{-C n/\log\log n}$.
\end{corollary}
\begin{proof}
We make use of a result of Gr\"unbaum and Moore \cite{GM} (see \S\ref{sec:previous}) which states
that if $n$ is odd, $E \subseteq \ZZ_n$, and $\ft{\chi_E}(1) \neq 0$ then
$E$ is determined by its $3$-deck. The rest follow from Corollary \ref{cor:at-one} with $k=1$.
\end{proof}

For arbitrary $n$ we lose a little in the exponent. Probably this is unnecessary.
\begin{theorem}\label{th:exceptions}
If $E$ is a random subset of $\ZZ_n$ the probability that $E$ is not determined
by its $3$-deck is at most
$$
2^{-C_\epsilon n^{1-\epsilon}},
$$
for any $\epsilon>0$.
\end{theorem}

For the proof we use some notions 
(recall Definition~\ref{extdom} and Notation~\ref{notation}) 
and lemmas from \S \ref{sec:positive} and also some new ones. 
Write
$$
A_x = \Set{k \in \ZZ_n: (k,n) \le x},
$$
and write ${\rm GAP}(B)$ for the size of the largest interval
contained in the complement of $B \subseteq \ZZ_n$.

\bl \label{GAP}
${\rm GAP}(A_{d(n)}) \le d(n)$, where $d(n)$ denotes
the number of divisors of $n$.
\el

\bp
Suppose that $I =\Set{a, a+1, \ldots, a+d(n)} \subseteq A_{d(n)}^c$
is an interval of size $d(n)+1$,
and $i, j \in I$, $i \neq j$.
Then $(i,n)>d(n)$ and $(j,n)>d(n)$.
It follows that $(i,n) \neq (j,n)$, otherwise we would have
$\Abs{i-j} \ge (i,n) > d(n)$, which cannot happen as all distances in $I$
are at most $d(n)$.
Thus, to each $i \in I$ there corresponds a different divisor of $n$, namely $(i,n)$.
But this cannot happen as $I$ has $d(n)+1$ members but there are only $d(n)$
different divisors of $n$.
\ep

\begin{lemma}\label{gapd}
If $\{0,1,\ldots,d\}\su A\su\Z_n$ and ${\rm GAP}(A)\le d$
then $A$ is an extendable domain.
\end{lemma}

\begin{proof}
Let $h:A\to\R/\Z$ be an additive function. 
First by induction we get that $h(j)=j h(1)\ (\md 1)$ for $j=1,\ldots,d$. 
Since ${\rm GAP}(A)\le d$, this can be extended by induction to all of $A$.
\end{proof}

\begin{lemma}\label{dninsupp}
If $E\su\Z_n$ and $\{1,2,\ldots,d(n)\}\su\sFTE$
then $E$ is determined up to translation by its $3$-deck.
\end{lemma}

\begin{proof}
By Fact~\ref{classes}, 
$\{1,2,\ldots,d(n)\}\su\sFTE$ implies that
$A_{d(n)}\su\sFTE$.
Thus by Lemma~\ref{GAP},
${\rm GAP}(\sFTE)\le d(n)$. Hence by Lemma~\ref{gapd}, $\sFTE$ is 
an extendable domain. Therefore by Lemma~\ref{domain}~(1),
$E$ is determined up to translation by its $3$-deck.
\end{proof}

\begin{proof}(Theorem~\ref{th:exceptions})
By Corollary \ref{cor:at-one},
\begin{eqnarray*}
\Prob{\exists j \in \Set{1,2,\ldots,d(n)}: \ft{\chi_E}(j) = 0}
 & \le & d(n) 2^{-C n/d(n)\log\log n}\\
 & \le & C_\epsilon n^\epsilon 2^{-C_\epsilon n^{1-\epsilon}}\\
 & \le & 2^{-C_{\epsilon'}n^{1-\epsilon'}},
\end{eqnarray*}
where $\epsilon'>0$ is again arbitrary,
and we used the fact that $d(n) = O(n^\epsilon)$ for all $\epsilon>0$ \cite{HW}.
By Lemma~\ref{dninsupp}, 
this completes the proof of Theorem~\ref{th:exceptions}.
\end{proof}

\end{document}